\documentclass[11pt, leqno]{amsart}
\usepackage{amsmath,amssymb,amsthm}
\usepackage{graphicx}
\usepackage{subfigure}
\usepackage{enumerate}
\usepackage{fullpage}
\usepackage{url}
\usepackage{float}
\usepackage{xspace}
\usepackage{hyperref}
\usepackage{color}
\usepackage{thmtools, thm-restate}
\usepackage{etoolbox}
\usepackage[page]{appendix}

\newtheorem{theorem}{Theorem}

\newtheorem{proposition}{Proposition}

\newtheorem{remark}{Remark}

\title{A short solution of the kissing number problem in dimension three}

\author{Alexey Glazyrin}
\address{School of Mathematical \& Statistical Sciences, The University of Texas Rio Grande Valley, Brownsville, TX 78520}
\email{alexey.glazyrin@utrgv.edu}

\begin{document}

\maketitle
\begin{abstract}

In this note, we give a short solution of the kissing number problem in dimension three.

\end{abstract}

\section{Introduction}

The problem of finding the maximum number of non-overlapping unit spheres tangent to a given unit sphere is known as \textit{the kissing number problem}. Sch{\"u}tte and van der Waerden \cite{sch53} settled the thirteen spheres problem (the kissing number problem for dimension three) that was the subject of the famous discussion between Isaac Newton and David Gregory in 1694. A sketch of an elegant proof was given by Leech \cite{lee56}. The thirteen spheres problem continues to be of interest to mathematicians, and new proofs have been published in recent years \cite{mae07, bor03b, ans04, mus06}. In other dimensions, the kissing number problem is solved only for $d=1, 2, 8, 24$ \cite{lev79,odl79}, and for $d=4$ \cite{mus08a}.

The proof presented in the note is similar to the proof in \cite{mus06} and the solution of the kissing problem in dimension four \cite{mus08a} (see also \cite{pfe07}) but the auxiliary function is chosen more carefully so the case analysis is much simpler and there is almost no spherical geometry required to cover the cases which seems to be the main component of all known proofs so far.

\begin{theorem}\cite{sch53}\label{thm:kiss}
The kissing number in dimension three is 12.
\end{theorem}

For our proof, we use ideas from the linear programming approach. The method was discovered by Delsarte \cite{del73} for the Hamming space, then extended to the spherical case \cite{del77} and generalized by Kabatyansky and Levenshtein \cite{kab78}. We use the properties of Gegenbauer polynomials defined recursively as follows.
$$G_0^{(d)}(t)=1,\ \ \ G_1^{(d)}(t)=t,\ \ \ G_k^{(d)}(t)=\frac {(d+2k-4)\, t\, G_{k-1}^{(d)}(t) - (k-1)\, G_{k-2}^{(d)}(t)}{d+k-3}.$$
In particular, the Delsarte method in the spherical case is based on the following proposition.
\begin{proposition}\cite{del77,kab78}\label{prop:delsarte}
For any finite set $X=\{x_1,\ldots, x_N\}\subset \mathbb{S}^{d-1}$ and any $k\geq 0$,
$$\sum\limits_{1\leq i,j \leq N}G_k^{(d)} (\langle x_i,x_j \rangle) \geq 0.$$
\end{proposition}

\section{A short proof of Theorem \ref{thm:kiss}}

\begin{figure}
\centering
\includegraphics[width=0.5\linewidth]{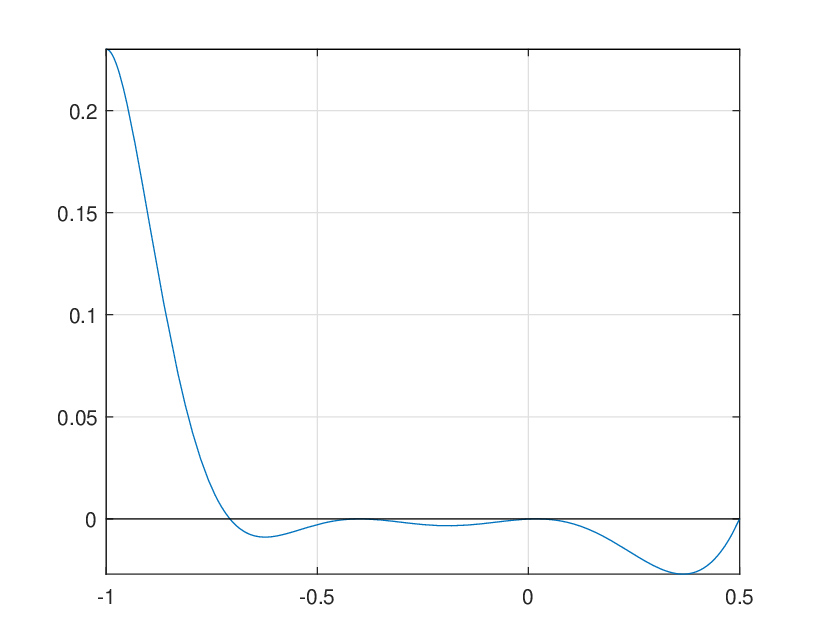}
\caption{Plot of $f(t)$ for $t\in [-1,1/2]$.}
\label{fig:graph}
\end{figure}

Let $f(t)=0.09465869 + 0.17273741\, G_1^{(3)}(t) + 0.33128438\, G_2^{(3)}(t) + 0.17275228\, G_3^{(3)}(t) +$\\ $0.18905584\, G_4^{(3)}(t) + 0.00334265\, G_5^{(3)}(t) + 0.03616728\, G_9^{(3)}(t)$ (see Figure \ref{fig:graph} for the plot of $f(t)$).

Assume we have $N$ non-overlapping unit spheres tangent to a given unit sphere $\mathbb{S}^2$. Then all pairwise angular distances between points of tangency $x_1, \ldots, x_N$ in $\mathbb{S}^2$ are at least $\pi/3$. If we show that for each $i$, $\sum_{j=1}^N f(\langle x_i,x_j \rangle) \leq 1.23$ then we can conclude the statement of the theorem. Indeed, on the one hand $\sum_{i,j=1}^N f(\langle x_i,x_j \rangle) \leq 1.23 N$. On the other hand, Proposition \ref{prop:delsarte} implies $\sum_{i,j=1}^N f(\langle x_i,x_j \rangle) \geq \sum_{i,j=1}^N 0.09465869 = 0.09465869 N^2$. Therefore, $N\leq 1.23/0.09465869 \approx 12.99405263$.

Fix $x=x_i$. The polynomial $f$ is negative on $[-1/\sqrt{2}, 1/2]$ so the positive contribution to the sum $\sum_{j=1}^N f(\langle x,x_j \rangle)$ can be made only by points $x_j$ in the open spherical cap $C$ with center $-x$ and angular radius $\pi/4$. No more than 3 points with pairwise angular distances at least $\pi/3$ can fit in $C$. Indeed, if there were at least 4 points $y_1, y_2, y_3, y_4$ in $C$ then at least one angle $\angle (y_i,-x, y_j)$ would be no greater than $\pi/2$. By the spherical law of cosines, the angular distance between $y_i$ and $y_j$ would be less than $\pi/3$.

If there is exactly one point $y$ in $C$, then
$$f(1)+f(\langle x,y \rangle) \leq f(1) + \max\limits_{t\in[-1,-1/\sqrt{2}]} f(t) \leq 1.23.$$

For two points $y,z$ in $C$, the angular distance between $y$ and $-x$ is at least $\frac \pi {12}$ by the triangle inequality for $y, z, -x$. Hence if $\langle x,y \rangle=t$ then $t$ cannot be less than $-\cos \frac \pi {12} $. By the triangle inequality, $\langle x,z \rangle\geq \alpha(t)=\frac 1 2 t-\frac {\sqrt{3}}2 \sqrt{1-t^2}$. Since $f$ is decreasing on $I=[-\cos \frac \pi {12},-1/\sqrt{2}]$,

$$f(1)+f(\langle x,y \rangle) + f(\langle x,z \rangle) \leq f(1) + \max_{t\in I} (f(t)+ f(\alpha(t)))\leq 1.23.$$

For three points $y,z,w$ in $C$, we use the monotonicity of $f$ on $I$ and move them as close as possible to $-x$. This way we get at least two of the three pairwise angular distances equal to $\pi/3$. Assume $\langle y,z \rangle=\langle z,w \rangle=1/2$. Note that $z,w,x$ cannot belong to the same great circle because otherwise $y$ does not fit in $C$. This means we can always  move $w$ keeping $\langle w,z \rangle=1/2$ and decreasing $\langle x,w \rangle$. The process stops in two possible cases: $w$ reaches the boundary of $C$ or $\langle y,w \rangle$ becomes 1/2. In the former case we are left with the case of two points in $C$ covered above. Now we can assume that $\langle y,z \rangle=\langle z,w \rangle=\langle y,w \rangle=1/2$.

We denote $t=\langle x,y \rangle, v=\langle x,z \rangle, u=\langle x,w \rangle$ for the rest of the proof. Without loss of generality, $t\leq v\leq u$. We note that $t\geq -2\sqrt{2}/3$. Otherwise, the circle with center $y$ and angular radius $\pi/3$ would intersect $C$ by an arc such that the distance between its endpoints is less than $\pi/3$. We also note that $t\leq -\sqrt{2/3}$ because otherwise $f(1)+f(t)+f(v)+f(u)\leq f(1)+3f(-\sqrt{2/3})\leq 1.23$.

For certain values of $u$, we can obtain the required inequality by bounding $f$ above by a linear function. Note that $t+v+u=\langle x, y+z+w\rangle\geq -||y+z+w|| = -\sqrt{6}$. For $t\in[-2\sqrt{2}/3,-0.76]$, $f(t)\leq -0.95 t - 0.6995$. Then, if we assume $u\leq -0.76$, 

$$f(1)+f(t)+f(v)+f(u)\leq f(1)-0.95(t+v+u)-2.0985\leq f(1)+0.95\sqrt{6}-2.0985\leq 1.23.$$

From now on, we assume $u\geq -0.76$ and, subsequently, $f(u)\leq f(-0.76)$. We fix $y$ and move $z$ and $w$ trying to minimize $v=\langle x,z \rangle$ and keeping the conditions $\langle y,z \rangle=\langle z,w \rangle=\langle y,w \rangle=1/2$, $u=\langle x,w\rangle \leq -1/\sqrt{2}$ (we disregard the condition $t\leq v\leq u$ for now). Since $z$ moves on a circle, this process may stop only in two cases: when $x, y, z$ belong to the same great circle or when $u$ becomes $-1/\sqrt{2}$. The former case is impossible as then $w$ would not fit in $C$. In the latter case, we can find the maximal possible value of $v$ explicitly by considering the Gram matrix of the vectors $x, y, z, w$ and equating its determinant to 0: $v=\beta(t)=\frac 1 3(t-1/\sqrt{2}-\sqrt{-8t^2-4\sqrt{2}t+2})$ (the second root of the quadratic equation gives values larger than $-1/\sqrt{2}$), where $t\in J = [-2\sqrt{2}/3,-\sqrt{2/3}]$. Then
$$f(1)+f(t)+f(v)+f(u) \leq f(1) + \max_{t\in J} (f(t)+ f (\beta(t))) + f(-0.76)\leq 1.23.$$

%

\begin{remark}
The function $f(t)$ was found by using a fixed value of 1.23 and maximizing the constant term in the Gegenbauer expansion while imposing several required conditions. In particular, the conditions $f(t)\leq 0$ on $[-1/\sqrt{2},1/2]$, $f(1)+f(t)\leq 1.23$ on $[-1,-1/\sqrt{2}]$, $f(1)+2f(-\sqrt{3}/2)\leq 1.23$, $f(1)=1$, $f(-1)=0.23$ were used. The optimization problem was then solved using SeDuMi optimization package in MATLAB. Coefficients were rounded to 8 decimal digits (after this the last two conditions stopped being true). All inequalities in the proof are easily verifiable. For convenience, their explicit forms are available in a separate file attached to the manuscript and at the author's website.
\end{remark}

\section{Acknowledgments}

The author thanks three anonymous referees whose help allowed him to fix the case of three points and improve the readability of the paper in general.

\bibliographystyle{amsplain}

\end{document}